# Finiteness Properties of Soluble $S$-Arithmetic Groups – a Survey

Kai-Uwe Bux

## 1 Introduction: Groups and Geometries

Every group is supposed to act upon a certain set preserving some additional structure. This set, which is almost always a space, should be associated to the group in a natural way – in many cases, the set comes first and the group is associated to the set. Linear groups act on vector spaces. Fuchsian groups act on the hyperbolic plane. Symmetry groups of geometric configurations act upon these objects.

In some cases the right space to act on is not that easily found. The mapping class group of a surface, i.e., the group of homotopy classes of homeomorphisms of that surface does not act upon the surface since it is a proper quotient of its automorphism group. Nevertheless, it acts upon the Teichmüller space of the surface. The group of outer automorphisms of a free group of finite rank acts on the Culler-Vogtmann space (outer space). These spaces have been christened in honor of those who found them. This indicates that it is not at all easy to find the right space.

How can we distinguish the right action from other actions? Most groups admit actions on many spaces. What distinguishes a good choice? As a rule of thumb, we will aim at small stabilizers as well as a small quotient space. In the case of a free action, every orbit is isomorphic to the group, and if the quotient space is a point, then there is only one orbit to deal with. Hence, at least philosophically, small quotient and stabilizers will force the space acted upon to look roughly like the group itself. Hence one may be able to deduce some properties of the group by proving some analogues for the space.

### 1.1 Finiteness Properties

To give a flavour of the theorems one obtains by employing geometrically arising group actions, we shall deal with finite generation and finite presentation first. Both of the following theorems are due to A.M. Macbeath [Macb64, Theorem 1]. J.-P. Serre also gave a proof [Serr77, Chapitre 1 § 3 Appendice, page 45].

Let $X$ be topological space and $U$ be an open subset of $X$. We assume that $X$ is locally connected and locally simply connected. Hence covering space theory applies.

Suppose that the group $G$ acts on $X$ such that

$$X = GU := \bigcup_{g \in G} gU.$$

Put $\mathcal{H} := \{h \in G \mid hU \cap U \neq \emptyset\}$. Note that $\mathcal{H}$ is symmetric, i.e., $h^{-1} \in \mathcal{H}$ for every $h \in \mathcal{H}$.

**Theorem 1.1.** *If $X$ is nonempty and path connected, $\mathcal{H}$ generates $G$. In particular, if $\mathcal{H}$ is finite, then $G$ is finitely generated.*

We can refine this result to obtain presentations as follows. Let $\mathcal{X} := \{x_h \mid h \in \mathcal{H}\}$ be a set of letters, one for each element of $\mathcal{H}$. Fix the following set of relations $\mathcal{R} := \left\{x_{h_1} x_{h_2} x_{h_3}^{-1} \,\middle|\, h_i \in \mathcal{H},\ h_1 h_2 = h_3\right\}$. The group

$$\tilde{G} := \langle \mathcal{X} \mid \mathcal{R} \rangle$$

admits an obvious homomorphism $\varphi : \tilde{G} \to G$ taking $x_h$ to $h$.

**Theorem 1.2.** *Suppose $U$ is path connected and $X$ is 1-connected, i.e., nonempty, connected, and simply connected. Then $\varphi$ is an isomorphism.*

**Corollary 1.3.** *If $\mathcal{H}$ is finite in addition to the assumptions of Theorem 1.2, then $G$ has a finite presentation.*

So we see that finite generation or finite presentability of a group are somehow geometric or topological properties of the group.

**Example 1.4.** The integers $\mathbb{Z}$ acts on $\mathbb{R}$ by translations. The translates of the open interval $(-\frac{1}{2}, 1\frac{1}{2})$ cover $\mathbb{R}$. We find $\mathcal{H} = \{-1, 0, 1\}$ and obtain the finite presentation

$$\mathbb{Z} = \langle x_{-1}, x_0, x_1 \mid x_{-1}x_0 = x_0x_{-1} = x_{-1},\ x_0x_0 = x_{-1}x_1 = x_1x_{-1} = x_0,\ x_1x_0 = x_0x_1 = x_1 \rangle .$$

**Example 1.5.** Let $X = \mathbb{E}^2$ be the Euclidean plane and consider the group $G$ of those isometries of $X$ that leave invariant the tiling by equilateral triangles shown below. $G$ is a Euclidean Coxeter group.

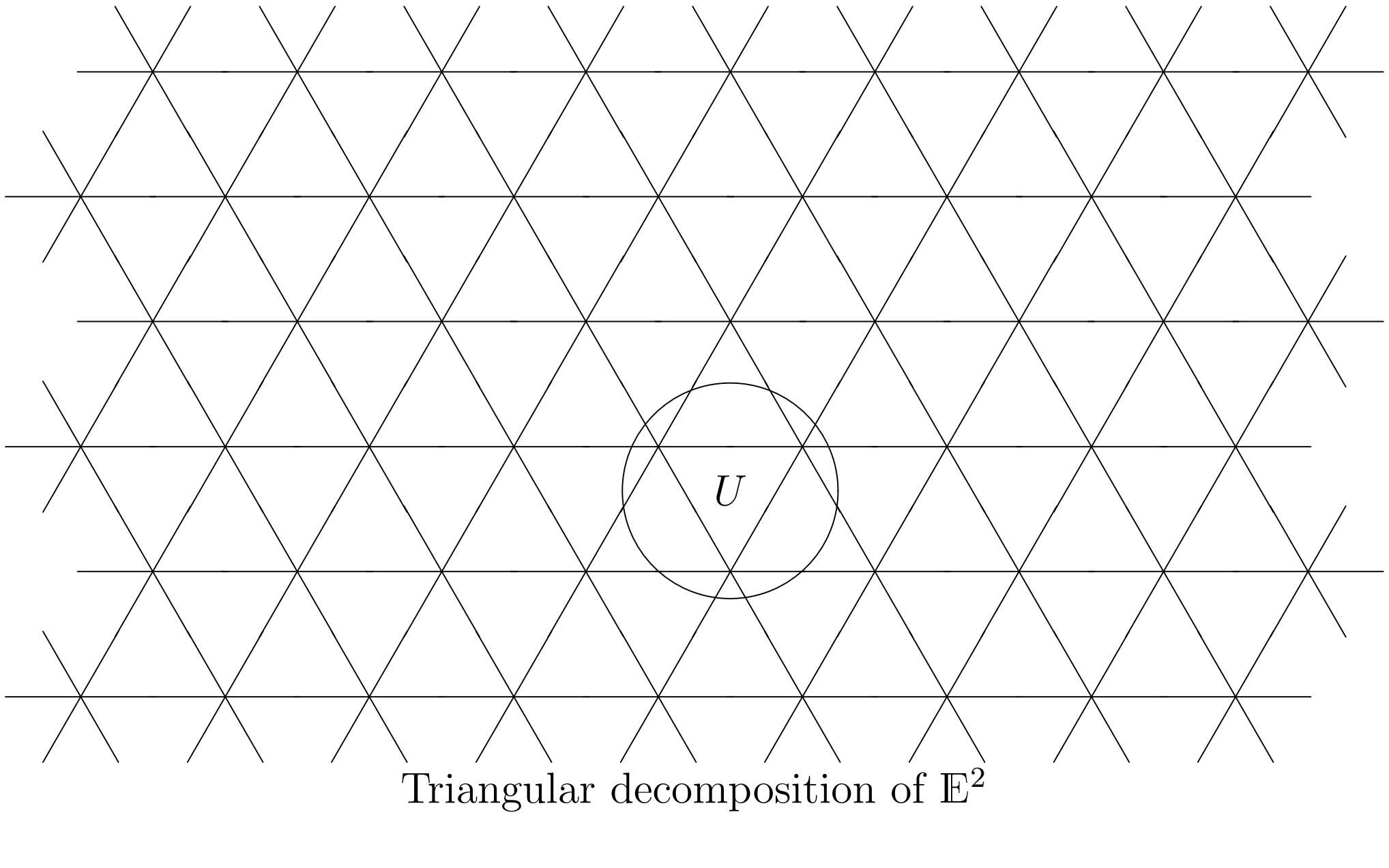

Triangular decomposition of $\mathbb{E}^2$

Let $U$ be an open disc containing a fixed triangle $\Delta$. Since an element of $G$ being an isometry is uniquely determined by the images of the three vertices of $\Delta$, only finitely many elements of $G$ can take $\Delta$ to a nearby triangle. Hence the assumptions of Theorem 1.2 are satisfied and we conclude that $G$ is finitely presented.

**Example 1.6.** More interesting is that we can do the very same thing for $\mathrm{SL}_2(\mathbb{Z})$. This group acts on the upper half plane $\mathbb{H}^2 = \{z \in \mathbb{C} \mid \Im(z) > 0\}$ via

$$\begin{pmatrix} a & b \\ c & d \end{pmatrix} z = \frac{az+b}{cz+d}.$$

The upper half plane is a model for hyperbolic geometry and the action of $\mathrm{SL}_2(\mathbb{Z})$ is by hyperbolic isometries. It leaves invariant a set of hyperbolic geodesics as shown in the figure – recall that, within the upper half plane model, geodesic lines are represented by vertical lines and half circles centered on the real line. Thus, we obtain a decomposition into *fundamental domains*.

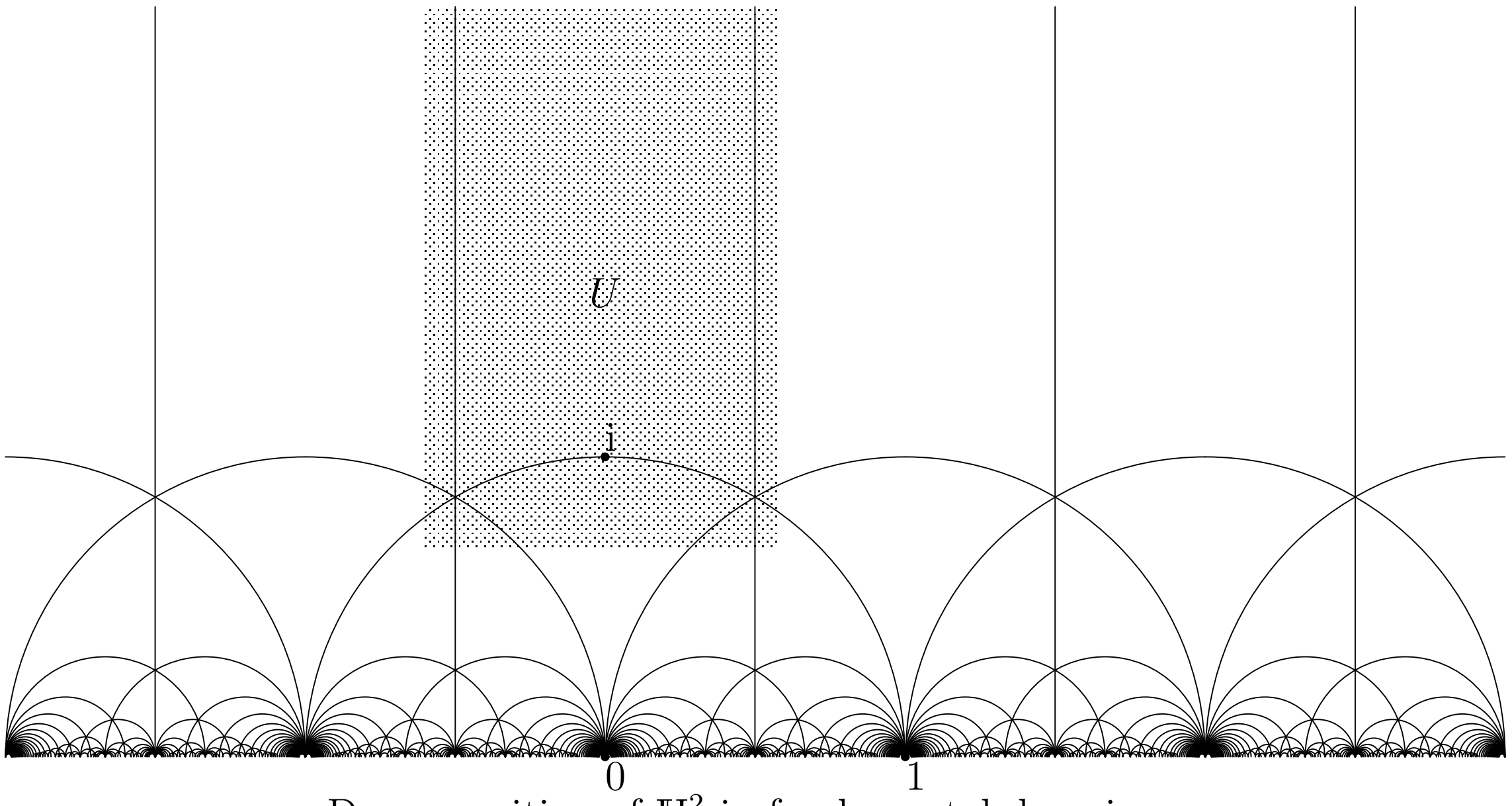

Decomposition of $\mathbb{H}^2$ in fundamental domains

Taking an open neighbourhood $U \subset \mathbb{H}^2$ of one of these domains and applying Theorem 1.2 shows that the group $\mathrm{SL}_2(\mathbb{Z})$ is finitely presented.

In fact even more is true. At every vertex of the decompositions there are three intersecting lines. So six geodesic rays issue from each vertex. Three of these run straight away to to infinity (that is, they approach the ideal boundary of the hyperbolic plane by going vertically upwards or approaching the real line) whereas the other three join the vertex to its three neighbouring vertices. The edges (line segments of finite length joining neighbouring vertices) form a three-valent tree.

The group $\mathrm{SL}_2(\mathbb{Z})$ acts on this tree since an isometry cannot take an edge to an infinite geodesic ray. Furthermore, $\mathrm{SL}_2(\mathbb{Z})$ has a subgroup of finite index that acts freely on this tree. Finally, the action of this subgroup has a compact quotient. This is one of the strongest finiteness conditions a group can satisfy. A group $G$ is *of finite type* if it acts cocompactly (i.e., with compact quotient) and freely on a contractible CW-complex. An equivalent condition is that $G$ is the fundamental group of a finite complex whose homotopy groups in dimensions $\geq 2$ vanish. $G$ is *virtually of finite type* if it has a subgroup of finite index which is of finite type.

As we have seen, $\mathrm{SL}_2(\mathbb{Z})$ is virtually of finite type. We will see below that this statement generalises to all arithmetic groups for which $\mathrm{SL}_2(\mathbb{Z})$ is the most simple nontrivial example.

There are finiteness properties of intermediate strength between finite generation and finite presentability on the one hand and being of finite type on the other hand. We can weaken the condition of being of finite type by restricting the cocompactness condition to a low dimensional skeleton: A group $G$ is *of type $F_m$* if it acts freely on a contractible CW complex whose $m$-skeleton has finite quotient mod $G$. The cellular chain complex of such a complex provides a free resolution of the trivial $\mathbb{Z}G$-module $\mathbb{Z}$ which is finitely generated up to dimension $m$. We can define a further weakening by saying that $G$ is *of type $FP_m$* if $\mathbb{Z}$ considered as a trivial $\mathbb{Z}G$-module has a projective resolution that is finitely generated up to dimension $m$. Hence every group of type $\mathrm{F}_m$ is of type $\mathrm{FP}_m$. The converse, however, does not hold [BeBr97, Examples 6.3(3)]. A group is finitely generated if and only if it is of type $\mathrm{F}_1$ (or $\mathrm{FP}_1$; the difference does not show up in this dimension), and it is finitely presented if and only if it is of type $\mathrm{F}_2$.

Recall that the $m$-skeleton of a contractible CW complex is $(m-1)$-*connected*, that is, its homotopy groups are trivial in dimensions less than $m$.

**Theorem 1.7 ([Brow87, Proposition 1.1]).** *Suppose $G$ acts cocompactly on an $(m-1)$-connected CW complex $X$ by cell permuting homeomorphisms such that the stabilizer of each cell $c$ is of type $F_{m-\dim(c)}$. Then $G$ is of type $F_m$.*

In particular, this theorem allows for finite cell stabilizers whereas the definition of higher finiteness properties uses a free action.

Although the relationship of Theorem 1.7 and Corollary 1.3 is not apparent, the former generalises the latter. Let $G$, $X$, and $U$ be as in Theorem 1.2. The *nerve* of the covering $X = \bigcup_{g\in G} gU$ is the simplicial complex $N$ defined by the following two conditions.

1. The vertex set of $N$ is $G$.

2. A finite subset $\sigma \subseteq G$ spans a simplex if and only if

$$U_\sigma := \bigcap_{g\in\sigma} gU \neq \emptyset.$$

The group $G$ acts on $N$ on the left by multiplication. The stabilizer of any simplex is finite. The action is cocompact if the set $\mathcal{H} = \{h \in G \mid U \cap hU \neq \emptyset\}$ is finite. Hence Corollary 1.3 follows from Theorem 1.7 and the following

**Proposition 1.8.** *If $X$ is 1-connected and $U$ is path connected then the geometric realization $|N|$ is 1-connected.*

**Proof.** Our proof is modeled on an argument of Quillen's [Quil78, Section 7].

We may assume that $X$ is locally path connected. Otherwise we retopologize $X$ by the topology generated by path components of open subsets of $X$. With this topology, $X$ is locally path connected and still 1-connected [Span66, Chapter 2, Exercises A2-3]. The subset $U$ remains open and path connected.

Let $X$ be a topological space and $\mathcal{X}$ be the category of open subsets in $X$ with inclusions as morphisms. Given a partially ordered set $D$, a *good cover of $X$ over $D$* is a functor $\mathcal{U} : D \to \mathcal{X}$ satisfying $X = \varinjlim_{\alpha\in D} \mathcal{U}(\alpha)$. A *$\mathcal{U}$-covering* of $X$ is a covering space of $X$ that covers every $\mathcal{U}(\alpha)$ evenly. Let $\mathcal{S}\mathrm{et}^{\mathrm{bij}}$ be the category of sets with bijections as the only morphisms. There is a one-to-one correspondence between functors $\boldsymbol{E} : D \to \mathcal{S}\mathrm{et}^{\mathrm{bij}}$ and $\mathcal{U}$-coverings $\pi : E \to X$ given by

$$E = \varinjlim_{\alpha\in D} \mathcal{U}(\alpha) \times \boldsymbol{E}(\alpha).$$

The covering $\pi : E \to X$ covers $X$ evenly if and only if $\boldsymbol{E}$ is naturally equivalent to a constant functor.

The set of simplices in the nerve $N$, also denoted by $N$, is partially ordered by the face relation. The functor $\boldsymbol{U}$ assigning to every simplex $\sigma \in N$ the set

$$\boldsymbol{U}(\sigma) := U_\sigma = \bigcap_{g\in\sigma} gU$$

is a good cover of $X$ over $N$. Since $X$ is 1-connected, every covering of $X$ covers $X$ evenly and is a $\boldsymbol{U}$-covering for that reason. Hence every functor $N \to \mathcal{S}\mathrm{et}^{\mathrm{bij}}$ is naturally equivalent to a constant functor.

There also is a good cover $\boldsymbol{N}$ of the geometric realization $|N|$ over $N$ given by

$$\boldsymbol{N}(\sigma) := \mathrm{St}^{\circ}_{N}(\sigma).$$

Since the open star of every simplex is contractible, any covering of $|N|$ covers every $\boldsymbol{N}(\sigma)$ evenly. Hence all coverings of $|N|$ are $\boldsymbol{N}$-coverings. These, corresponding to functors from $N$ to $\mathcal{S}\mathrm{et}^{\mathrm{bij}}$ and therefore to coverings of $X$, cover $|N|$ evenly. Hence $|N|$ is 1-connected. **q.e.d.**

In Theorem 3.1 we will quote a refinement of Theorem 1.7 that provides necessary and sufficient conditions for higher finiteness properties.

There are many other finiteness properties. For instance, a group could be of finite cohomological or geometrical dimension. See [Brow82, Chapter VIII] for more about this.

**Example 1.9.** Finally, let us consider a group that leads us directly to the heart of the matter. Fix a finite field $k$ and let $K := k(t)$ be the field of rational functions over $k$. This field contains the ring $L := k[t,t^{-1}]$ as a subring. We will consider the group

$$\mathrm{B}_2^0(L) := \left\{ \left( \begin{array}{cc} u & p \\ 0 & u^{-1} \end{array} \right) \middle| \, u \in L^*, \; p \in L \right\} \subseteq \mathrm{SL}_2(L).$$

It acts on $L$ via

$$\left( \begin{array}{cc} u & p \\ 0 & u^{-1} \end{array} \right) q := \frac{uq+p}{u^{-1}} = u^2 q + up$$

as inspired by the action of $\mathrm{SL}_2(\mathbb{Z})$ on $\mathbb{C}$ in the preceding example, and as above, we will turn this action into an action on a tree.

For each $m \in \mathbb{Z}$ put

$$V_m := \left\{ \sum_{i \in \mathbb{Z}} \alpha_i t^i \in L \,\middle|\, \alpha_i = 0 \; \forall i \leq m \right\}$$

We regard these sets as neighbourhoods of $0 \in L$. Translation yields neighbourhoods for each Laurent polynomial $p \in L$

$$V_m(p) := \{ q \in L \,|\, p - q \in V_m \}.$$

We call $m$ the *radius* of $V_m(p)$. Let

$$\mathcal{V} := \{ V_m(p) \,|\, p \in L, \; m \in \mathbb{Z} \}$$

be the set of all these neighbourhoods partially ordered by inclusion. Observe that $\mathrm{B}_2^0(L)$ acts on $\mathcal{V}$ preserving the ordering.

The inclusion relation gives rise to a directed graph with vertex set $\mathcal{V}$. For two elements $V, V' \in \mathcal{V}$, insert an edge pointing from $V$ to $V'$ if there is no $W \in \mathcal{V}$ with $V' \subsetneqq W \subsetneqq V$ while $V' \subsetneqq V$. So the edges point from $V$ precisely to those elements of $\mathcal{V}$ that are maximal with respect to inclusion in $V$. Note that these are finite in number whence we obtain a locally finite graph $T_+$, on which $\mathrm{B}_2^0(L)$ acts.

$T_+$ does not contain directed cycles, and given two neighbourhoods $V$, $V' \in \mathcal{V}$ with nonempty intersection $V \cap V' \neq \emptyset$, one of them contains the other. Hence each vertex of $T_+$ is the terminal vertex of at most one edge, though $|k|$ edges issue from there. It follows that $T_+$ does not contain undirected cycles. Furthermore, every neighbourhood $V$ is contained in some neighbourhood $V_m$ of 0. Hence every vertex in $T_+$ is joined by an edge path to the line

$$\cdots \rightarrow V_{m-1} \rightarrow V_m \rightarrow V_{m+1} \rightarrow \cdots$$

whence $T_+$ is a tree.

Since the action of $\mathrm{B}_2^0(L)$ on $T_+$ does not have finite stabilizers, we cannot obtain finite generation of this group by applying Theorem 1.2 to $T_+$. Nevertheless, there is a remedy. The "co-neighbourhoods"

$$W_m := \left\{ \sum_{i \in \mathbb{Z}} \alpha_i t^i \in L \,\middle|\, \alpha_i = 0 \; \forall i \geq m \right\}$$

give rise to a second tree $T_-$ in a completely analogous manner. The diagonal action of $\mathrm{B}_2^0(L)$ on the product $T_+ \times T_-$ permutes the cells (vertices, edges, and squares) of this complex. Stabilizers of cells are now finite, however, the action has no compact quotient. So we pass to an appropriate subspace.

As neighbourhoods have radii, co-neighbourhoods have "co-radii". The map that associates to every vertex $(V_m(p), W_n(q))$ in in the Product $T_+ \times T_-$ the number $m-n$ extends linearly to edges and squares of $T_+ \times T_-$. We denote this map by

$$\boldsymbol{\pi} : T_+ \times T_- \to \mathbb{R}$$

and observe that it is invariant under the action of $\mathrm{B}_2^0(L)$. The space

$$X := \boldsymbol{\pi}^{-1}(0)$$

is connected, and Theorem 1.7 implies that $\mathrm{B}_2^0(L)$ is finitely generated.

This construction will be considered again in Section 4 where we shall also prove that $\mathrm{B}_2^0(L)$ is not finitely presented.

## 1.2 Arithmetic Groups

We will think of a *linear algebraic group* as a group of matrices of determinant 1 which is defined by polynomial equations that the matrix coefficients are supposed to satisfy. These coefficients can be taken from any fixed ring containing the constants that occur in the defining equations. A typical example is $\mathrm{SO}_r$. Note that a set of polynomial equations chosen at random will almost always fail to define a subgroup of the special linear group SL. Nevertheless, if it does, then the fact that the determinant is 1 shows that if the coefficients of a matrix belong to a certain ring, the coefficients of its inverse will do so, too. [Bore91] may serve as a reference on linear algebraic groups.

Among the linear algebraic groups, there are two extreme types: soluble groups on the one hand side and semi-simple groups on the other hand side. A linear algebraic group is *semi-simple* if it does not contain a nontrivial connected soluble normal subgroup. Therefore, any linear algebraic group is an extension of a semi-simple group by a soluble group. Reductive groups are close to semi-simple groups. A linear algebraic group is *reductive* if it contains only "small" connected soluble normal subgroups, that is, these soluble normal subgroups do not contain unipotent elements. Maximal connected soluble subgroups of linear algebraic group are called *Borel subgroups*. Since they need not be normal, semi-simple and reductive groups may contain many Borel subgroups. They are a most important tool in studying reductive groups.

An *arithmetic group* or an *$S$-arithmetic* is obtained from a linear algebraic group when the coefficients of its matrices are chosen from an *arithmetic ring* or an *$S$-arithmetic ring*. We already have encountered such rings in the examples. $\mathbb{Z}$ is the most simple arithmetic ring. Generally, any ring of algebraic integers in a number field is called arithmetic. In the number field case, e.g., over the rationals $\mathbb{Q}$, $S$-arithmetic rings are obtained from arithmetic rings by localizing at a finite set $S$ of prime elements. That is, one allows some primes to be inverted. All these rings live inside of *global number fields* which are by definition finite extensions of the rationals $\mathbb{Q}$.

Number fields have a nice arithmetic structure, and it has been recognized long ago that there are fields of positive characteristic which exhibit a very similar behaviour. These are the so called *global function fields* which are finite extensions of fields $k(t)$ of rational functions over finite fields $k$. It does not make sense to look for algebraic integers in these fields – they form a finite subfield, the so called *constant functions*. Hence there are no interesting arithmetic rings. Nevertheless, $S$-arithmetic subrings can be defined even though the notion of primes has to be recasted in terms of valuations. The precise definitions, which apply to number fields as well, are given in the next section. At the moment it suffices to know that the polynomial ring $k[t]$ and the ring $k[t, t^{-1}]$ of Laurent polynomials are $S$-arithmetic rings over a set of one prime and two primes, respectively.

$\mathrm{SL}_2(\mathbb{Z})$ and $\mathrm{B}_2^0(k[t, t^{-1}])$ are typical examples. $\mathrm{SL}_2(\mathbb{Z})$ is an arithmetic subgroup of a semi-simple linear algebraic group and $\mathrm{B}_2^0(k[t, t^{-1}])$ is an $S$-arithmetic subgroup of a Borel group.

As long as people have studied arithmetic groups, they have been investigating finiteness properties. This is not surprising since, as we have seen, finiteness properties are of geometric nature and the geometries associated to algebraic groups are a fundamental tool in studying them and their arithmetic subgroups. This goes for $S$-arithmetic subgroups as well. For an arithmetic group over a number field, the symmetric space of the real Lie group of real points of the algebraic group is a geometric model whereas, over function fields, Bruhat-Tits buildings provide good geometries for $S$-arithmetic groups. The hyperbolic plane of Example 1.6 is a symmetric space, and the trees $T_+$ and $T_-$ of Example 1.9 are Bruhat-Tits buildings associated to the two primes defining the ring of Laurent polynomials.

Let us discuss the number field case first. Let $\mathcal{G}$ be a linear group defined over the global number field $K$ with ring of integers $\mathcal{O}$. Beside finite generation and finite presentability (see [BoHa62] and [Behr62]), the first result on higher finiteness properties is a theorem of M.S. Raghunatan [Ragh68, Theorem 1, Corollaries 2 and 4] implying that an arithmetic subgroup of $\mathcal{G}(K)$ is of type $\mathrm{F}_\infty$ provided that $\mathcal{G}$ is semi-simple. He proved even more, namely that any arithmetic subgroup of a semi-simple group is virtually of finite type. This is a far-reaching generalisation of Example 1.6.

A. Borel and J.-P. Serre reproved Raghunathan's result by different means [BoSe73, Theorem 9.3]. They also observed that this theorem implies that an arithmetic subgroup of $\mathcal{G}(K)$ is of type $\mathrm{F}_\infty$ regardless of whether the linear group $\mathcal{G}$ is semi-simple. Furthermore, they generalised the result to $S$-arithmetic groups [BoSe76, Proposition 6.10]. In this case however, they had to assume that the linear algebraic group is reductive. Then, as above, an $S$-arithmetic subgroup is virtually of finite type.

What about non-reductive linear algebraic groups? M. Kneser proved that $\mathcal{G}(\mathcal{O}_S)$ is finitely presented if and only if for each non-Archimedian prime $v \in S$ the locally compact group $\mathcal{G}(K_v)$ is compactly presented where $K_v$ denotes the completion of $K$ at $v$ [Knes64]. Compact presentability is the analogue of finite presentability in the category of locally compact topological groups. Using Kneser's result, H. Abels characterized all finitely presented $S$-arithmetic groups over number fields [Abel87]. The main point is that $\mathcal{G}(K_v)$ is compactly presented if and only if a maximal $K_v$-split soluble subgroup $\mathcal{B}(K_v)$ is compactly presented. This reduces the problem to soluble groups. For these, he proved:

**Theorem 1.10.** *Let $\Gamma$ be a soluble $S$-arithmetic group. Then there is a short exact sequence*

$$N \hookrightarrow \Gamma \twoheadrightarrow Q$$

*where $N$ is nilpotent and $Q$ contains a finitely generated Abelian subgroup of finite index. $\Gamma$ is finitely presented if and only if the following two conditions hold*

1. *The Abelianization ${}^{N}\!/_{[N,N]}$ is a tame $\mathbb{Z}Q$-module.*
2. $\mathrm{H}_2(N;\mathbb{Z})$ *is finitely generated over $\mathbb{Z}\Gamma$.*

To explain tameness, we have to consider homomorphisms $\chi : Q \to \mathbb{R}$. Such a homomorphism defines a submonoid $Q_\chi := \{q \in Q \mid \chi(q) \geq 0\}$. A $\mathbb{Z}Q$-module is said to be *tame* if, for every homomorphism $\chi$, the module is finitely generated either over $\mathbb{Z}Q_\chi$ or over $\mathbb{Z}Q_{-\chi}$ or over both rings. This concept was introduced by R. Bieri and R. Strebel in [BiSt80] to study finite presentability of metabelian groups.

This line of thought has been generalised by H. Abels and A. Tiemeyer to deal with higher finiteness properties. It yields the following

**Theorem A.** *Let $K$ be a global number field and $S$ a finite nonempty set of primes containing all Archimedian primes. Let $\mathcal{O}_S$ be the corresponding $S$-arithmetic subring of $K$. Furthermore let $\mathcal{B}$ be a Borel subgroup scheme of a reductive group defined over $K$. Then the $S$-arithmetic group $\mathcal{B}(\mathcal{O}_S)$ is of type $F_\infty$.*

P. Abramenko [unpublished] proved this result before. His proof establishes that $\mathcal{B}(\mathcal{O}_S)$ is virtually of finite type. In spite of these results, there does not yet exist a generalisation of Theorem 1.10, i.e., a list of necessary and sufficient conditions for a soluble $S$-arithmetic group to be of type $\mathrm{F}_m$. Further examples of soluble arithmetic groups are due to H. Abels (see [Abel79] and [AbBr87]). We will discuss the proof of Theorem A Section 3 following Abels and Tiemeyer.

Less is known in the function field case. Even for reductive groups the problem of determining the finiteness properties of $S$-arithmetic subgroups is unsolved. For these, H. Behr has given a complete solution for finite generation and finite presentability in [Behr98]. Concerning higher finiteness properties, there are mainly two series of examples. U. Stuhler proved that $\mathrm{SL}_2(\mathcal{O}_S)$ is of type $\mathrm{F}_{|S|-1}$ but not of type $\mathrm{FP}_{|S|}$ ([Stuh76] and [Stuh80]) where $\mathcal{O}_S$ is an $S$-arithmetic subring of any global function field. This series shows a positive influence of the number of primes. On the other hand, P. Abramenko [Abra87] and H. Abels [Abel91], have proved independently, that $\mathrm{SL}_n(\mathbb{F}_q[t])$ is of type $\mathrm{F}_{n-2}$ but not of type $\mathrm{FP}_{n-1}$ provided that $q$ is big enough. P. Abramenko generalised this series to other classical groups [Abra96]. It is not yet known whether the assumption on $q$ can be dropped. H. Behr, though, has proposed a strategy for proving a positive answer recently [Behr99].

For soluble groups, the author obtained the following result ([Bux97a] and [Bux97b], see [Bux99] for an English version).

**Theorem B.** *Let $K$ be a global function field and $S$ a finite nonempty set of primes. Let $\mathcal{O}_S$ be the corresponding $S$-arithmetic subring of $K$. Furthermore let $\mathcal{B}$ be a Borel subgroup scheme of a Chevalley group defined over $\mathbb{Z}$. Then the $S$-arithmetic group $\mathcal{B}(\mathcal{O}_S)$ of type $F_{|S|-1}$ but not of type $FP_{|S|}$.*

Concerning other soluble groups, little is known in the function field case (see [Bux97b, Bemerkungen 8.10 and 8.11] or [Bux99, Remarks 8.6 and 8.7]). We will outline the proof of Theorem B in Section 4. Since the detailed proofs will be published elsewhere, we focus on the geometric aspects.

# 2 Notation

Let $K$ be a *global field*. Hence $K$ is either a *global number field* (that is, a finite extension of the rationals $\mathbb{Q}$) or a *global function field* (that is a finite extension of the field of rational functions $k(t)$ of a finite field $k$). For detailed information about global fields and their arithmetic, see [O'Me73], [CaFr67], or [Weil73]. Let $A$ be the finite set of non-Archimedian primes of $K$, which is empty if $K$ is a function field. Recall that a *prime* or *place* $v$ of $K$ is an equivalence class of valuations. To this we can associate the corresponding completion $K_v$ of $K$. This is a local field. Its additive group therefore has a Haar measure. For each $x \in K_v$, the multiplication by $x$ induces an automorphism of the additive group, which rescales the Haar measure by a positive real number $|x|_v$. If $v$ is non-Archimedian, the corresponding *valuation ring* $\mathcal{O}_v := \{x \in K_v \mid |x|_v \leq 1\}$ of $v$-adic integers is a local ring. For any finite, nonempty set $S$ of primes containing $A$ let $\mathcal{O}_S := \{f \in K \mid |f|_v \leq 1 \ \forall v \notin S\}$ denote the corresponding $S$-arithmetic ring and $\mathbb{A}_S := \times_{v \in S} K_v \times \times_{v \notin S} \mathcal{O}_v$ the ring of $S$-adeles. The adele ring is the direct limit $\mathbb{A} := \varinjlim_S \mathbb{A}_S$. This is a locally compact ring. It contains $K$ as a discrete subring via the diagonal embedding because each element in $K$ belongs to $\mathcal{O}_v$ for all but finitely many primes $v$ – in the language of the function field case: a function has at most finitely many poles. Let $\mathcal{G}$ be a linear algebraic group defined over $K$ represented as a group of matrices with determinant 1. Then for any subring $R$ of any $K$-algebra, the group $\mathcal{G}(R)$ of $R$-points is defined. Hence we regard $\mathcal{G}$ as a functor taking (topological) subrings of $K$-algebras to (topological) groups. The group $\mathcal{G}(\mathcal{O}_S)$ is called an $S$-arithmetic subgroup of $\mathcal{G}$. It does depend on the matrix representation chosen for $\mathcal{G}$. Nevertheless, any two representations of $\mathcal{G}$ yield commensurable $S$-arithmetic groups. For this reason, any group commensurable with $\mathcal{G}(\mathcal{O}_S)$ is called *$S$-arithmetic*, as well. In the number field case, the set $A$ of Archimedian primes is nonempty. Hence, over number fields, $A$-arithmetic groups exist. These are simply called *arithmetic*.

Since commensurable groups have identical finiteness properties, we may confine our investigation to groups of the form $\mathcal{G}(\mathcal{O}_S)$.

# 3 The Number Field Case

Throughout this section, $K$ is a number field. H. Abels and A. Tiemeyer introduced higher compactness properties $\mathrm{C}_m$ and $\mathrm{CP}_m$ for locally compact groups which extend compact generation (equivalent to $\mathrm{C}_1$ and to $\mathrm{CP}_1$) and compact presentability (equivalent to $\mathrm{C}_2$) [AbTi97]. For a discrete group, the properties $\mathrm{C}_m$ (respectively $\mathrm{CP}_m$) and $\mathrm{F}_m$ (respectively $\mathrm{FP}_m$) are equivalent. In [Tiem97], a Hasse Principle is derived relating finiteness properties of $S$-arithmetic groups to compactness properties of linear algebraic groups over local fields. We will discuss this in Section 3.2. Applying this, Tiemeyer settled Theorem A mainly using the technique of contracting automorphisms already present in [Abel87]. This approach yields a comparatively simple proof, also due to Tiemeyer, for the result of Borel´s and Serre´s on the finiteness properties of $S$-arithmetic groups. We will present this argument in Section 3.3. It does not yield, however, that these groups are virtually of finite type.

## 3.1 Finiteness Properties and Compactness Properties

The Starting point of the definition of compactness properties is the following celebrated criterion of K. Brown's.

**Theorem 3.1 ([Brow87, Theorems 2.2 and 3.2]).** *Let $\Gamma$ be a discrete group and $(X_\alpha)_{\alpha\in D}$ a directed system of cocompact $\Gamma$-CW-complexes with equivariant connecting morphisms and $(m-1)$-connected limit. Suppose that, for each $X_\alpha$, the stabilizers of all $j$-cells are of type $F_{m-j}$ (respectively $FP_{m-j}$).*

*Then $\Gamma$ is of type $F_m$ (respectively $FP_m$) if and only if the directed systems of homotopy groups $(\pi_n(X_\alpha))_{\alpha\in D}$ (respectively the directed systems of reduced homology groups $(\mathrm{H}_n(X_\alpha))_{\alpha\in D}$) are essentially trivial for all $n < m$.*

Herein, a directed system $(G_\alpha)_{\alpha\in D}$ of groups is called *essentially trivial* if for each element $\alpha \in D$ there is a $\beta \geq \alpha$ such that the natural map $G_\alpha \to G_\beta$ is trivial.

Defining

$$\mathrm{Mor}\left((G_\alpha)_{\alpha\in D}, (G'_\beta)_{\beta\in D'}\right) := \varprojlim_{\alpha\in D} \varinjlim_{\beta\in D'} \mathrm{Mor}\left(G_\alpha, G'_\beta\right)$$

we turn the class of all directed systems of groups into the category $\mathrm{ind}_{\mathrm{Gr}}$ [AbTi97, Section 1]. In $\mathrm{ind}_{\mathrm{Gr}}$, the essentially trivial directed systems of groups are precisely those that are isomorphic to the initial element of this category represented by the trivial system of trivial groups. Henceforth we consider directed systems of groups as elements of $\mathrm{ind}_{\mathrm{Gr}}$.

For every group $\Gamma$, there is a canonical directed system of cocompact simplicial $\Gamma$-complexes indexed by the finite subsets $F \subseteq \Gamma$, namely the $\Gamma$-orbit of $F$ in the the simplicial complex of all finite subsets of $\Gamma$, acted upon by left translation. Simplex stabilizers are finite and the limit is contractible. Hence finiteness properties can – in principle – be tested with this directed system in view of Theorem 3.1. We will somehow mimic this construction for locally compact groups.

Let $G$ be a locally compact group, $X$ a locally compact space together with an action of $G$ on $X$, and $\mathrm{E}X$ the free simplicial set over $X$ whose $m$ skeleton is just $X^{m+1}$. The face (respectively degeneration) maps are given by deleting (respectively doubling) a coordinate. Note that $G$ acts on $\mathrm{E}X$ diagonally and therefore also on its geometric realization $|\mathrm{E}X|$. This is a contractible space. Furthermore, let $\mathcal{C}X$ denote the set of all compact subsets of $X$ directed by inclusion. For each set $C \in \mathcal{C}X$, $G \cdot \mathrm{E}C := \{(gx_0, \ldots, gx_m) \mid g \in G,\ \{x_0, \ldots, x_m\} \subseteq C,\ m \in \mathbb{N}\}$ is a $G$-invariant simplicial subset of $\mathrm{E}X$. Hence we have a directed system $\mathrm{E}_G X := (G \cdot \mathrm{E}C)_{C \in \mathcal{C}X}$ of simplicial $G$-sets depending functorially on $G$ and $X$, whose limit is $\mathrm{E}X$.

Taking the geometric realization first and passing afterwards to either homotopy groups $\pi_n$ or reduced homology $\mathrm{H}_n$, induces the directed systems of groups $\mathrm{ind}{-}\pi_n(\mathrm{E}_G X) := (\pi_n(|G \cdot \mathrm{E}C|))_{C \in \mathcal{C}X} \in \mathrm{ind}_{\mathrm{Gr}}$ and $\mathrm{ind}{-}\mathrm{H}_n(\mathrm{E}_G X) := (\mathrm{H}_n(|G \cdot \mathrm{E}C|))_{C \in \mathcal{C}X} \in \mathrm{ind}_{\mathrm{Gr}}$, respectively.

**Definition 3.2.** The pair $(G, X)$ satisfies *condition $P_m$* (respectively *condition $PP_m$*) if $\mathrm{ind}{-}\pi_n(\mathrm{E}_G X) \cong 0$ (respectively $\mathrm{ind}{-}\mathrm{H}_n(\mathrm{E}_G X) \cong 0$) for all $n < m$. That is, these associated directed systems are essentially trivial.

Now we can define compactness properties of locally compact groups. The group $G$ is *of type $C_m$* (respectively *of type $CP_m$*) if the pair $(G, G)$ satisfies condition $\mathrm{P}_m$ (respectively $\mathrm{PP}_m$) where $G$ is considered as a $G$-space via left multiplication.

It is immediate from Brown's criterion that these notions agree for discrete groups with the properties $\mathrm{F}_m$ and $\mathrm{FP}_m$, respectively.

**Remark 3.3.** Recall that an object $R$ in a category is a retract of the object $O$ if there are arrows $R \hookrightarrow O$ and $O \twoheadrightarrow R$ whose composition is the identity on $R$. Retract diagrams are preserved by covariant and contravariant functors. Hence a retract of the group $G$ enjoys at least the same compactness and finiteness properties as $G$ because of the categorial characterisation – vanishing of some functorially assigned objects – of compactness and finiteness properties given above. **q.e.d.**

We will use the following lemmata.

**Lemma 3.4 ([AbTi97, Lemma 3.2.2]).** *The $P_m$- and $PP_m$-conditions are independent of the choice of the $G$-space $X$ provided the action of $G$ on $X$ is proper. Hence the conditions $C_m$ and $CP_m$ can be tested using any proper $G$-space.*

**Corollary 3.5 ([AbTi97, Corollary 3.2.3]).** *If $G$ is a locally compact group and $B \leq G$ is a closed subgroup with compact quotient ${}^{G}/_{B}$ then $G$ and $B$ have the same compactness properties.*

In the case of finiteness properties, one can use Brown's Criterion 3.1 to slightly improve the lemma.

**Lemma 3.6.** *Let $X$ be a $G$-space such that, for all $k$-tuples of points in $X$, the intersection of their stabilizers is of type $F_{m-k}$. Then $G$ is of type $F_m$ (respectively $FP_m$) if and only if the pair $(G, X)$ satisfies $P_m$ (respectively $PP_m$).* **q.e.d.**

**Lemma 3.7 ([AbTi97, Lemma 3.1.1]).** *A finite direct product of locally compact groups is of type $C_m$ (respectively $CP_m$) if each of its factors is of type $C_m$ (respectively $CP_m$).*

## 3.2 The Hasse Principle

In this section we shall outline the proof of

**Theorem 3.8 (Hasse Principle [Tiem97, Theorem 3.1]).** *Let $\mathcal{G}$ be a linear algebraic group. Then the $S$-arithmetic group $\mathcal{G}(\mathcal{O}_S)$ is of type $F_m$ (respectively $FP_m$) if and only if for each non-Archimedian prime $v \in S$ the locally compact completion $\mathcal{G}(K_v)$ is of type $C_m$ (respectively $CP_m$).*

**Proof.** For a non-Archimedian prime $v$, the group $\mathcal{G}(\mathcal{O}_v)$ is a compact open subgroup of $\mathcal{G}(K_v)$, which is a locally compact group.

Let us put $\mathcal{G}_S := \times_{v\in S}\,\mathcal{G}(K_v)$. This is a locally compact group since $S$ is finite. The subgroup $K_S := \times_{v\in S-A}\,\mathcal{G}(\mathcal{O}_v)$ is compact and open in $\mathcal{G}_{S-A}$.

Since we have

$$\mathcal{G}(\mathbb{A}_S) = \bigtimes_{v\in S}\mathcal{G}(K_v) \times \bigtimes_{v\notin S}\mathcal{G}(\mathcal{O}_v)$$

we can, on the one hand, consider $K_S$ as a subgroup of $\mathcal{G}(\mathbb{A}_S)$; on the other hand, there is a natural projection map $\pi_S : \mathcal{G}(\mathbb{A}_S) \to \mathcal{G}_S$.

Since $K_S$ is open in $\mathcal{G}_{S-A}$ $X := {}^{\mathcal{G}_{S-A}}\!/_{K_S}$ is a discrete set upon which $\mathcal{G}_{S-A}$ acts properly. $\mathcal{G}(\mathcal{O}_S)$ also acts on $X$, via $\pi_S$. The stabilizers are of the form $\mathcal{G}(\mathcal{O}_S) \cap gK_Sg^{-1}$ and can be shown to be commensurable with $\mathcal{G}(\mathcal{O}_S) \cap K_S = \mathcal{G}(A)$ whence they are arithmetic and therefore of type $\mathrm{F}_\infty$. Hence, by Lemma 3.6, $\mathcal{G}(\mathcal{O}_S)$ is of type $\mathrm{F}_m$ (respectively $\mathrm{FP}_m$) if and only if the pair $(\mathcal{G}(\mathcal{O}_S), X)$ is of type $\mathrm{P}_m$ (respectively $\mathrm{PP}_m$).

On the other hand, by Lemma 3.7, $\mathcal{G}(K_v)$ is of type $\mathrm{C}_m$ (respectively $\mathrm{CP}_m$) for each $v \in S - A$ if and only if the direct product $\mathcal{G}_{S-A} = \times_{v\in S-A}\,\mathcal{G}(K_v)$ is of type $\mathrm{C}_m$ (respectively $\mathrm{CP}_m$). Because of Lemma 3.4, this in turn holds if and only if the pair $(\mathcal{G}_{S-A}, X)$ is of type $\mathrm{P}_m$ (respectively $\mathrm{PP}_m$).

Hence we have to compare the actions of $\mathcal{G}(\mathcal{O}_S)$ and $\mathcal{G}_{S-A}$ on $X$. In particular it suffices to show that the filtrations $\mathrm{E}_{\mathcal{G}(\mathcal{O}_S)}X$ and $\mathrm{E}_{\mathcal{G}_{S-A}}X$ are cofinal. So let $K$ be a compact (that is finite) subset of $X$. Obviously

$$\mathcal{G}(\mathcal{O}_S)\cdot \mathrm{E}K \le \mathcal{G}_{S-A}\cdot \mathrm{E}K$$

whence we only have to prove that there is another finite subset $L \subset X$ such that

$$\mathcal{G}_{S-A}\cdot \mathrm{E}L \le \mathcal{G}(\mathcal{O}_S)\cdot \mathrm{E}K.$$

In [Tiem97, Corollary 2.3], it is deduced from Borel's Finiteness Theorem [Bore63, Theorem 5.1] that there are finitely many $g_1, \ldots, g_m$ such that $\mathcal{G}_{S-A} = \bigcup_{i=1}^m \mathcal{G}'(\mathcal{O}_S)g_iK_S$ where $\mathcal{G}'(\mathcal{O}_S) := \pi_S(\mathcal{G}(\mathcal{O}_S))$.

Observe that for finite $K \subset X$ the set $L := \bigcup_{i=1}^m g_iK_SK$ is finite since $K_S$ acts with finite orbits on $X$ as it is a compact group. Then, for every $g \in \mathcal{G}_{S-A}$, there is a $g_i$ and an element $\gamma \in \mathcal{G}(\mathcal{O}_S)$, acting via its image in $\mathcal{G}'(\mathcal{O}_S)$, such that $g \in \gamma g_iK_S$ whence $gK = \gamma g_iK_SK \subseteq \gamma L$. From this,

$$\mathcal{G}_{S-A}\cdot \mathrm{E}L \le \mathcal{G}(\mathcal{O}_S)\cdot \mathrm{E}K$$

follows immediately. **q.e.d.**

**Corollary 3.9.** *Let $\mathcal{P}$ be a parabolic subgroup of $\mathcal{G}$. Then $\mathcal{G}(\mathcal{O}_S)$ and $\mathcal{P}(\mathcal{O}_S)$ have the same finiteness properties.*

**Proof.** For each prime $v$, since ${}^{\mathcal{G}(K_v)}\!/_{\mathcal{P}(K_v)}$ is compact by [BoTi65, Proposition 9.3], the groups $\mathcal{G}(K_v)$ and $\mathcal{P}(K_v)$ have the same compactness properties by Corollary 3.5. Now the claim follows by the Hasse principle. **q.e.d.**

**Remark 3.10.** Corollary 3.9 yields a quick and dirty proof for Theorem A, since $S$-arithmetic subgroups of Chevalley groups are of type $\mathrm{F}_\infty$ by [BoSe76, Proposition 6.10].

Using Corollary 3.9 this way, however, is somehow approaching the problem from the wrong angle. The corollary should be considered the other way round, namely, as reducing the problem of determining the finiteness properties of general $S$-arithmetic groups to $S$-arithmetic subgroups of connected soluble linear algebraic groups.

## 3.3 Stand Alone Proofs

In this section we shall indicate how $S$-arithmetic subgroups of soluble groups can be studied. Let $G = T \ltimes U$ a locally compact semi-direct product of locally compact groups where $T$ is Abelian. An element $t \in T$ is *contracting* or acts *by contraction* on $U$ if its positive powers converge uniformly to the identity on compact subsets of $U$. The following proposition shows that one is fortunate having a contracting element at hand.

**Proposition 3.11 ([Tiem97, Theorem 4.3]).** *If $T$ contains a contracting element, $G$ and $T$ have identical compactness properties.*

This lemma generalises [Abel87, Proposition 1.3.1].

**Proof.** Since $T$ is a retract of $G$, we only have to show, that compactness properties of $T$ are inherited by $G$. So assume that $T$ is of type $\mathrm{C}_n$.

Let $K_G$ be a compact subset of $G$. Enlarging this subset if necessary, we might assume that is of the form
$$K_G = K_T K_U$$
where $K_T$ is compact in $T$ and $K_U$ is a compact neighbourhood of the identity in $U$ satisfying $tK_Ut^{-1} \subseteq K_U$ after replacing $t$ by a power if necessary. Since we can test compactness properties of $T$ by means of its action on $G$, there is a compact subset $L \supseteq K_G \cup K_G t^{-1}$ such that the inclusion of $T \cdot \mathrm{E}K_G \hookrightarrow T \cdot \mathrm{E}L$ induces trivial maps in homotopy groups up to dimension $n$. We claim that this holds also for the inclusion of $G \cdot \mathrm{E}K_G \hookrightarrow G \cdot \mathrm{E}L$.

Consider the map
$$\begin{aligned} \alpha_t : \mathrm{E}G &\rightarrow \mathrm{E}G \\ (g_0, \ldots, g_k) &\mapsto \left(g_0 t^{-1}, \ldots, g_k t^{-1}\right) \end{aligned}$$
and observe that it takes $G \cdot \mathrm{E}K_G$ to itself because any simplex therein is of the form
$$\sigma = (gx_0u_0, \ldots, gx_ku_k)$$

with $x_i \in K_T$ and $u_i \in K_U$ and is taken to

$$\begin{aligned}\alpha_t(\sigma) &= \left(gx_0u_0t^{-1}, \ldots, gx_ku_kt^{-1}\right)\\ &= \left(gx_0t^{-1}tu_0t^{-1}, \ldots, gx_kt^{-1}tu_kt^{-1}\right)\\ &= \left(gt^{-1}x_0tu_0t^{-1}, \ldots, gt^{-1}x_ktu_kt^{-1}\right)\end{aligned}$$

which belongs to $G \cdot \mathrm{E}K_G$ because of $tu_it^{-1} \in K_U$. From the first line of these equations, moreover, it follows that $\sigma$ and $\alpha_t(\sigma)$ are faces of a common simplex in $G \cdot \mathrm{E}L$ since $K_Gt^{-1} \subseteq L$. Hence a sphere $\$$ in $G \cdot \mathrm{E}K_G$ is homotopic to the sphere $\alpha_t(\$)$ inside $G \cdot \mathrm{E}L$.

The map $\alpha_t$, however, has another remarkable property. It contracts $\mathrm{E}G$ towards $T \cdot \mathrm{E}K_G$: Consider a vertex $xu$ in $\mathrm{E}G$, $x \in T$, $u \in U$. We have

$$\alpha_t{}^n(xu) = t^{-n}xt^nut^{-n}$$

and eventually $t^nut^{-n} \in K_U$ since $t$ is contracting. Hence, any compact subset of $\mathrm{E}G$ is taken to $T \cdot \mathrm{E}K_G$ by some positive power of $\alpha_t$.

Combining the properties of $\alpha_t$, we can move any sphere $\$$ of dimension $\leq n$ in $G \cdot \mathrm{E}K_G$ inside $G \cdot \mathrm{E}L$ into $T \cdot \mathrm{E}K_G$. This sphere, however, is homotopically trivial inside $T \cdot \mathrm{E}L \subseteq G \cdot \mathrm{E}L$ by choice of $L$. **q.e.d.**

**Proposition 3.12.** *Let $\mathcal{P}$ be a parabolic subgroup of the reductive group $\mathcal{G}$, both defined over a local field $F$ of characteristic $0$ whose associated valuation $v : F^* \to \mathbb{R}$ is non-Archimedian. Then $\mathcal{P}(F)$ is of type $C_\infty$.*

**Remark 3.13.** Combined with the Hasse principle 3.8, this yields a proof of Theorem A which does not depend on [BoSe76]. Moreover, Corollary 3.9 implies that $\mathcal{G}(F)$ is of type $\mathrm{C}_\infty$, too.

**Proof.** In view of [BoTi65, Proposition 9.3] and Corollary 3.5, it suffices to construct a maximal connected $F$-split soluble subgroup $\mathcal{B} \leq \mathcal{G}$ such that $\mathcal{B}(F)$ is of type $\mathrm{C}_\infty$.

We start with a maximal $F$-split torus $\mathcal{T}$ in $\mathcal{G}$. The set of all characters $\mathcal{T} \to \mathfrak{Mult}$ is a finitely generated Abelian group, denoted by $X^+(\mathcal{T})$. Therefore $X^+(\mathcal{T}) \otimes_{\mathbb{Z}} \mathbb{R}$ is a real vector space. This is where the *root system* $\Phi := \Phi(\mathcal{G}, \mathcal{T})$ lives. It is the set of weights of the adjoint representation of $\mathcal{T}$ over the Lie algebra $\mathfrak{g}$ of $\mathcal{G}$. Thus, for every root $\alpha \in \Phi$, we are given a corresponding weight space and a corresponding unipotent subgroup $\mathcal{U}_\alpha \leq \mathcal{G}$. Fix a base of $\Phi$ thereby determining the positive half $\Phi^+$, and let $\mathcal{U}$ be the group generated by all $\mathcal{U}_\alpha$ for $\alpha \in \Phi^+$. This is a unipotent group [BoTi65, 3.8 (iv)], upon which $\mathcal{T}$ acts by conjugation inside $\mathcal{G}$ because we started with the adjoint representation. Furthermore, $\mathcal{B} := \mathcal{T} \ltimes \mathcal{U}$ is a maximal connected $F$-split soluble subgroup of $\mathcal{G}$. Hence, it suffices to show, that $\mathcal{T}(F) \ltimes \mathcal{U}(F)$ is of type $\mathrm{C}_\infty$ for each prime $v$.

To prove this, we will find a contracting element $t \in \mathcal{T}(F)$. Then the claim follows from Proposition 3.11 since $\mathcal{T}(F)$ is of type $\mathrm{C}_\infty$ for it contains a cocompact finitely generated Abelian group.

Whether an element of $\mathcal{T}(F)$ is contracting can be checked via its action on the Lie algebra: An element $t$ acts by contraction on $\mathcal{U}(F)$ if for every positive root $\alpha \in \Phi^+$, the inequality $v(\alpha(t)) > 0$ holds. Consider now the one parameter subgroups in $\mathcal{T}$. They form a lattice of maximum rank in the dual of $X^+(\mathcal{T}) \otimes_{\mathbb{Z}} \mathbb{R} \supset \Phi$. Hence there is a one parameter subgroup $\xi : \mathfrak{Mult} \to \mathcal{T}$ representing a point on which all positive roots are strictly positive. Finally, $t := \xi(x)$ fits our needs provided $v(x) > 0$. **q.e.d.**

# 4 The Function Field Case

From now on, $K$ is a global function field. In this case even $S$-arithmetic subgroups of reductive groups have finite finiteness length, and no Hasse Principle holds – the different primes have to work together in order to ensure higher finiteness properties. The main tool for establishing finiteness properties in this setting is to study the action of the $S$-arithmetic group $\mathcal{G}(\mathcal{O}_S)$ on the product of Bruhat-Tits buildings associated to $\mathcal{G}$ and the primes in $S$. Such buildings always exist if $\mathcal{G}$ is reductive [BrTi72]. For general $\mathcal{G}$, no replacement is known. Nevertheless, since Theorem B only deals with Borel subgroups of Chevalley groups, we can use the buildings associated to the latter in order to investigate the $S$-arithmetic subgroups of the former.

We shall freely use the terminology related to buildings and the reader should have a basic knowledge of these geometric objects and their relationship to algebraic groups. Standard references are [Brow89] and [Rona89]. Nevertheless, the main geometric argument in Section 4.2 does only depend on properties of buildings that are explicitly stated there.

## 4.1 Preliminaries on Chevalley Groups and their Associated Bruhat-Tits Buildings

Let $\mathcal{G}$ be a *Chevalley group*, i.e., a semi-simple linear algebraic group scheme defined over $\mathbb{Z}$. They have been introduced by Chevalley in [Chev60]. A standard reference is [Stei68]. Fix a Borel subgroup $\mathcal{B}$ of $\mathcal{G}$ and a maximal torus $\mathcal{T}$ inside $\mathcal{B}$ such that both are also defined over $\mathbb{Z}$.

Since a global function field $K$ has only non-Archimedean primes, there is a Euclidean Bruhat-Tits building $X_v := X(\mathcal{G}, K_v)$ associated to each prime $v \in S$. This is a CAT(0)-space whose boundary at infinity is the spherical building $\tilde{X}_v$ canonically associated to $\mathcal{G}$ and $K_v$. The group $\mathcal{B}(K_v)$ is the stabilizer of the *fundamental chamber* $C_v \in \tilde{X}_v$ at infinity. Likewise, $\mathcal{T}(K_v)$ stabilizes of the *standard apartment* $\Sigma_v$ in $X_v$, which corresponds to an apartment in $\tilde{X}_v$ containing $C_v$. It is a Euclidean Space in which we choose a simplicial cone $S_v \subseteq \Sigma_v$ representing $C_v$. We take the cone point to be the origin of $\Sigma_v$, turning it into a Euclidean vector space.

Recall that the root system $\Phi$ of $\mathcal{G}$ with respect to $\mathcal{T}$ consists of morphisms from $\mathcal{T}$ to $\mathfrak{Mult}$, the group scheme associating to each ring its group of multiplicative units, e.g., represented as diagonal $2 \times 2$ matrixes of determinant 1. So we can represent the root system on $\Sigma_v$ by a set of linear forms $\Phi_v$. Furthermore, we are given a base and a system of positive roots $\Delta \subseteq \Phi^+ \subset \Phi$ corresponding to the Borel subgroup $\mathcal{B}$.

These give rise to sets of linear forms $\Delta_v \subseteq \Phi_v^+ \subset \Phi_v$. The linear forms in $\Delta_v$ form a system of coordinates $\xi_v := \left(\alpha_v^{(1)}, \ldots, \alpha_v^{(r)}\right) : \Sigma_v \to \mathbb{R}^r$ on $\Sigma_v$ such that

$$S_v = \left\{ s_v \in \Sigma_v \,\middle|\, \alpha_v^{(j)}(s_v) \geq 0 \;\forall j \in \{1, \ldots, r\} \right\}$$

holds. This description remains true if all positive roots are taken into account.

We normalize the coordinates such that an element $t_v \in \mathcal{T}(K_v)$ acts on $\Sigma_v$ as a translation with coordinates

$$\left(\log(|\alpha^{(1)}(t_v)|_v), \ldots, \log(|\alpha^{(r)}(t_v)|_v)\right).$$

The subgroup of elements in $\mathcal{B}(K_v)$ that fix at least one point in $\Sigma_v$ is the group of $K_v$-points $\mathcal{U}(K_v)$ of the *unipotent radical* $\mathcal{U} := \ker(\mathcal{B} \to \mathcal{T})$. Moreover, $\Sigma_v$ is a strong fundamental domain for the action of $\mathcal{U}(K_v)$ on $X_v$, and there is a well defined map $\rho_v : X_v \to \Sigma_v$ called the *retraction centered at* $C_v$.

We will study the diagonal action of $\mathcal{B}(\mathcal{O}_S)$ on the product $\boldsymbol{X} := \times_{v \in S} X_v$. The retractions $\rho_v$ induce a retraction $\boldsymbol{\rho} : \boldsymbol{X} \to \boldsymbol{\Sigma} := \times_{v \in S} \Sigma_v$ onto the product of standard apartments.

The coordinates $\xi_v$ defined above give rise to a system of coordinates $\boldsymbol{\xi} : \boldsymbol{s} = (s_v)_{v \in S} \mapsto \left(\alpha_v^{(1)}, \ldots, \alpha_v^{(r)}\right)_{v \in S} \in \mathbb{R}^{r|S|}$ on $\boldsymbol{\Sigma}$. Since we normalized the coordinates $\xi_v$, the map $\boldsymbol{\zeta} : \boldsymbol{\Sigma} \to \mathbb{R}^r$, $\boldsymbol{s} = (s_v)_{v \in S} \mapsto \left(\sum_{v \in S} \alpha_v^{(1)}(s_v), \ldots, \sum_{v \in S} \alpha_v^{(r)}(s_v)\right)$ is invariant with respect to the action of $\mathcal{T}(\mathcal{O}_S)$ since the product formula [O'Me73, Theorem 33.1] implies that for every $x \in \mathcal{O}_S$,

$$\prod_{v \in S} |x|_v = 1.$$

Hence, the map $\boldsymbol{\pi} := \boldsymbol{\zeta} \circ \boldsymbol{\rho} : \boldsymbol{X} \to \mathbb{R}^r$ is invariant under the action of $\mathcal{B}(\mathcal{O}_S)$ on $\boldsymbol{X}$.

**Lemma 4.1 ([Bux97b, Lemma 2.3]).** *For every compact subset $C \subset \mathbb{R}^r$, the $S$-arithmetic group $\mathcal{B}(\mathcal{O}_S)$ acts cocompactly on $\boldsymbol{\pi}^{-1}(C)$.*

**Sketch of Proof.** We use adele topology. The main ingredient is that $\mathcal{U}(K)$ is a discrete subgroup of $\mathcal{U}(\mathbb{A})$ with compact quotient [Bux97b, Lemma 1.1]. This generalises the well known result that the additive group of $K$ is discrete in $\mathbb{A}$ and that $\mathbb{A}/K$ is compact [Weil73, Theorem 2, page 64]. It follows that for every polysimplex $\boldsymbol{\sigma} \in \boldsymbol{X}$ the double quotient

$$\mathcal{U}(\mathcal{O}_S) \backslash \mathcal{U}(\mathbb{A}_S) / \mathrm{Stab}_{\mathcal{U}(\mathbb{A}_S)}(\boldsymbol{\sigma})$$

is discrete and compact; hence finite. This implies that the map

$$\mathcal{U}(\mathcal{O}_S) \backslash \boldsymbol{X} \to \boldsymbol{\Sigma}$$

induced by $\rho$ is proper. Now the claim follows from Dirichlet's Unit Theorem [CaFr67, page 72]. **q.e.d.**

## 4.2 Positive Results

We outline the proof of the first half of Theorem B, namely

**Theorem 4.2.** *The group $\mathcal{B}(\mathcal{O}_S)$ is of type $F_{|S|-1}$.*

We apply Theorem 1.7. Lemma 4.1 exhibits subspaces of $\boldsymbol{X}$ on which $\mathcal{B}(\mathcal{O}_S)$ acts cocompactly. Cell stabilizers are finite as they are intersections of the compact stabilizers in $\mathcal{B}(\mathbb{A}_S)$ the discrete group $\mathcal{B}(\mathcal{O}_S)$. Hence it suffices to show that the preimage $\boldsymbol{\pi}^{-1}(C)$ is $(|S|-2)$-connected for some compact set $C \subset \mathbb{R}^r$, e.g., for a point. We put $\boldsymbol{Y} := \boldsymbol{\pi}^{-1}(0) = \boldsymbol{\rho}^{-1}(H) \subseteq \boldsymbol{X}$ where $H = \ker(\boldsymbol{\zeta} : \boldsymbol{\Sigma} \to \mathbb{R}^r)$.

We write $\boldsymbol{Y}$ as an ascending union

$$\boldsymbol{Y} = \bigcup_{j=0}^{\infty} \boldsymbol{Y}_j$$

where $\boldsymbol{Y}_0$ is contractible and $\boldsymbol{Y}_{j+1}$ is obtained from $\boldsymbol{Y}_j$ by gluing in a convex Euclidean set along its boundary such that, homotopically, this amounts to attaching a cell of dimension at least $|S|-1$. This way, no nontrivial homotopy elements are introduced up to dimension $|S|-2$.

### 4.2.1 Moufang Buildings and Λ-Complexes

The Moufang condition describes a certain interplay between the geometry of a building and its group of automorphisms. We shall not quote the precise definition here since Proposition 4.3 below, which is the technical core of the proof of Theorem 4.2, spells out all properties of Moufang buildings we use.

It is known that the buildings $X_v$ are locally finite Euclidean Moufang buildings [Bux97b, Fact 6.1]. Hence, the following proposition applies.

**Proposition 4.3 ([Bux97b, Lemma 6.2]).** *Let $X$ be a locally finite Euclidean Moufang building with a distinguished chamber $C$ at infinity and a distinguished apartment $\Sigma$ containing $C$. Then there is a sequence $\Sigma = \Sigma_0, \Sigma_1, \ldots$ of apartments such that the following conditions are satisfied:*

1. *Each $\Sigma_j$ contains $C$.*
2. *The $\Sigma_j$ cover $X$, i.e.,*
$$X = \bigcup_{j=0}^{\infty} \Sigma_j.$$
3. *For every $j > 0$, the new part*
$$N_j := \Sigma_j - \bigcup_{i<j} \Sigma_i$$
*is an intersection of open half apartments in $\Sigma_j$.*

*None of the $N_j$ contains $C$ because this chamber is already present in $\Sigma_0 = \Sigma$.*

**Remark 4.4.** The Moufang assumption is crucial to the proof of Proposition 4.3 given in [Bux97b], and I do not know whether the statement holds for non-Moufang locally finite Euclidean buildings. Since the buildings arising in the number field case are non-Moufang, a proof of Theorem A along the lines of the proof of Theorem B presented here does not yet exist.

We mention that A. von Heydebreck [Heyd99, Lemma 3.3] proved an exact analogue of Proposition 4.3 for finite buildings by purely geometric means that do not make use of any Moufang type assumption. For finite buildings, the distinguished chamber at infinity is to be replaced by a chamber inside the building.

Proposition 4.3 says that we can build $X$ starting from $\Sigma$ by attaching first the closure $\overline{N_1}$ along its boundary which lies in $\Sigma$. Then we glue in $\overline{N_2}$, $\overline{N_3}$, etc. Since $\overline{N_j}$ is an intersection of closed half apartments all of which avoid $C$, it can be retracted to its boundary. Hence the homotopy type does not change during this procedure. This recovers the important fact that the building $X$ is contractible.

Note that the retraction $\rho := \rho_{C,\Sigma} : X \to \Sigma$ restricts to an isomorphism of Coxeter complexes $\rho : \Sigma_j \to \Sigma$ for every $j$. Since these are Euclidean Coxeter complexes, this map is an isometry.

Since $\Sigma$ is a Euclidean Coxeter complex, the orthogonal part of its automorphism group is a finite reflection group associated to a root system $\Phi := \Phi(\Sigma)$, which we represent as a finite set of linear forms on $\Sigma$. The chamber $C$ at infinity determines a base and a subset of positive roots $\Phi^+$ such that $C$ is represented by the cone

$$\left\{ s \in \Sigma \,\middle|\, \alpha(s) \geq 0 \ \forall \alpha \in \Phi^+ \right\}.$$

We state some geometric axioms describing this setting.

**Definition 4.5.** Let $Y$ be a metric, piecewise Euclidean CW complex and $\pi : Y \to \mathbb{E}$ a projection from $Y$ onto a Euclidean vector space $\mathbb{E}$. Furthermore, let $\Lambda$ be a finite set of linear forms on $\mathbb{E}$. The pair $(Y, \pi)$ is called a $\Lambda$-*complex* if there exists a sequence $E_0$, $E_1$, ... of subcomplexes in $Y$ satisfying the following conditions:

1. The map $\pi$ restricts to an isometry $\pi|_{E_j} : E_j \to \mathbb{E}$ for every $j$.

2. The subcomplexes $E_j$ cover $Y$, i.e.,
$$Y = \bigcup_{j \geq 0} E_j.$$

3. For every $j > 0$, the new part
$$N_j := E_j - \bigcup_{i<j} E_i$$
is the interior of a subcomplex in $E_j$ whose $\pi$-image in $\mathbb{E}$ is of the form
$$\{e \in \mathbb{E} \mid \lambda(e) \leq c_\lambda\}$$
where $c_\lambda \in \mathbb{R} \cup \{\infty\}$ are constants not all of which are $\infty$. If $c_\lambda = \infty$, it imposes no restriction on $e$.

A sequence of this form is called *increasing*.

**Example 4.6.** Proposition 4.3 implies immediately that, keeping the notation from above, $(X, \rho : X \to \Sigma)$ is a $\Phi^+$-complex.

The following constructions provide many other examples. First, we consider direct products in order to deal with $\boldsymbol{X} = \times_{v\in S} X_v$.

**Lemma 4.7.** *For $i \in \{1,2\}$, let $(Y^i, \pi^i : Y^i \to \mathbb{E}^i)$ be an $\Lambda^i$-complex. Let $\mathrm{pr}_i : \mathbb{E}^1 \times \mathbb{E}^2 \to \mathbb{E}^i$ denote the canonical projection and put $\Lambda^1 \uplus \Lambda^2 := \{\lambda^i \circ \mathrm{pr}_i : \mathbb{E}^1 \times \mathbb{E}^2 \to \mathbb{R} \mid i \in \{1,2\}\,,\ \lambda^i \in \Lambda^i\}$. Then $(Y^1 \times Y^2, \pi^1 \times \pi^2 : Y^1 \times Y^2 \to \mathbb{E}^1 \times \mathbb{E}^2)$ is a $\Lambda^1 \uplus \Lambda^2$-complex.*

**Proof.** For $i \in \{1,2\}$, let $E_0^i, E_1^i, E_2^i, \ldots$ be an increasing sequences for $Y^i$. It is easy to check that $E_0^1 \times E_0^2$, $E_1^1 \times E_0^2$, $E_0^1 \times E_1^2$, $E_1^1 \times E_1^2$, $E_2^1 \times E_0^2$, $E_0^1 \times E_2^2$, $E_2^1 \times E_1^2$, $E_1^1 \times E_2^2$, $E_2^1 \times E_2^2$, $\ldots$ is an increasing sequence for $Y^1 \times Y^2$. **q.e.d.**

**Corollary 4.8.** *$(\boldsymbol{X} = \times_{v\in S} X_v, \boldsymbol{\rho} : \boldsymbol{X} \to \boldsymbol{\Sigma})$ is a $\biguplus_{v\in S} \Phi_v^+$-complex.*

**Proof.** All $X_v$ are locally finite Euclidean Moufang buildings [Bux97b, Fact 6.1]. **q.e.d.**

Now we turn to subcomplexes because we are interested in $\boldsymbol{Y} \subseteq \boldsymbol{X}$.

**Observation 4.9.** *Let $(Y, \pi : Y \to \mathbb{E})$ be a $\Lambda$-complex and let $\mathbb{E}'$ be a linear subspace of $\mathbb{E}$. Put $\Lambda|_{\mathbb{E}'} := \{\lambda|_{\mathbb{E}'} \mid \lambda \in \Lambda\}$.*

*Then $(\pi^{-1}(\mathbb{E}'), \pi\big|_{\pi^{-1}(\mathbb{E}')} : \pi^{-1}(\mathbb{E}') \to \mathbb{E}')$ is a $\Lambda|_{\mathbb{E}'}$-complex. We see this by taking an increasing sequence for $Y$ and intersecting it with the preimage of $\mathbb{E}'$.* **q.e.d.**

**Corollary 4.10.** *$(\boldsymbol{Y}, \boldsymbol{\rho}\big|_{\boldsymbol{Y}} : \boldsymbol{Y} \to H)$ is a $\left(\biguplus_{v\in S} \Phi_v^+\right)\big|_H$-complex.* **q.e.d.**

### 4.2.2 Connectivity of Λ-Complexes

We call a set $\Lambda$ of linear forms on a real vector *$m$-tame* if no positive linear combination of up to $m$ elements of $\Lambda$ vanishes, i.e., if 0 is not contained in the convex hull of up to $m$ elements of $\Lambda$.

**Proposition 4.11 ([Bux97a, Lemma 7.3]).** *Let $(Y, \pi : Y \to \mathbb{E})$ be a $\Lambda$-complex and suppose $\Lambda$ is $m$-tame. Then $Y$ is $(m-1)$-connected.*

**Proof.** Let $E_0, E_1, E_2, \ldots$ be an increasing sequence for $Y$. Put

$$Y_j := \bigcup_{i\leq j} E_j.$$

Since $Y$ is covered by the $E_j$, it suffices to show that every $Y_j$ is $(m-1)$-connected. Since $Y_0 = E_0$ is contractible, we can prove the claim by induction if we control the process of obtaining $Y_{j+1}$ from $Y_j$.

From

$$Y_{j+1} = Y_j \cup (E_{j+1} - Y_j)$$

and the fact that $E_0, E_1, E_2, \dots$ is increasing it follows that $Y_{j+1}$ is obtained from $Y_j$ by attaching, along its boundary, a convex set $N'_j$ of the form

$$N'_j = \{e \in \mathbb{E} \mid \lambda(e) \leq \lambda \ \forall \lambda \in \Lambda'\}$$

where $\Lambda' \subseteq \Lambda$ is a nonempty subset and $\lambda \in \mathbb{R}$.

We have to control the homotopy type of the pair $(N'_j, \partial(N'_j))$. If $N'_j = \partial(N'_j)$ there is nothing to do. Hence we can assume that $N'_j$ has nonempty interior.

First suppose that $\Lambda'$ spans the dual $\mathbb{E}^*$, that is,

$$0 = \bigcap_{\lambda \in \Lambda'} \ker(\lambda)\,.$$

There are two cases.

<u>$N'_j$ is unbounded</u>: The interior $N_j$ contains an infinite ray $r$, but it does not contain a whole line, since we assumed $0 = \bigcap_{\lambda \in \Lambda'} \ker(\lambda)$. Hence we may assume that $r$ starts at $x \in \partial(N'_j)$.

If the boundary $\partial(N'_j)$ does not contain a ray parallel to $r$ we can argue that moving all points parallel to $r$ towards the boundary $\partial(N'_j)$ defines a deformation retraction of $N'_j$ onto its boundary.

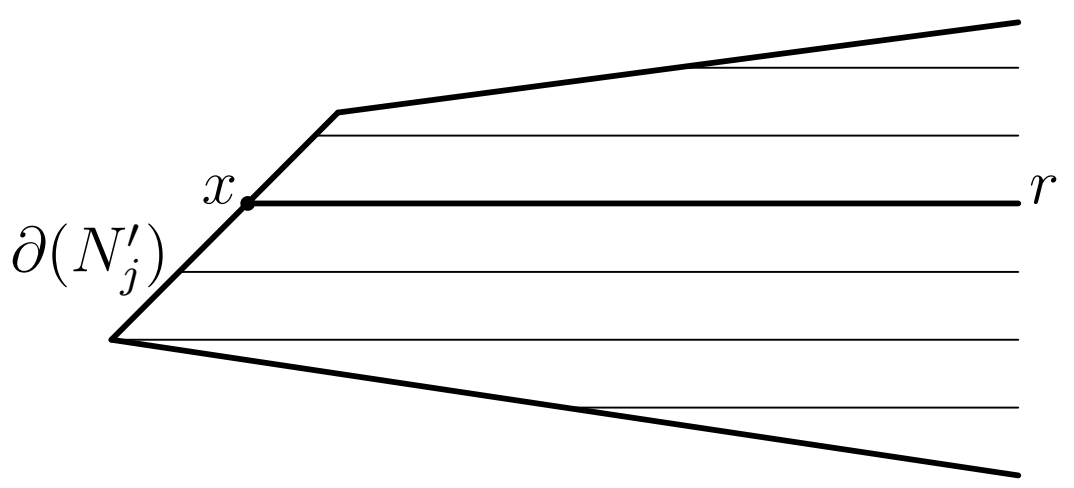

Good case

Unfortunately, this does not work if the boundary does contain a ray parallel to $r$.

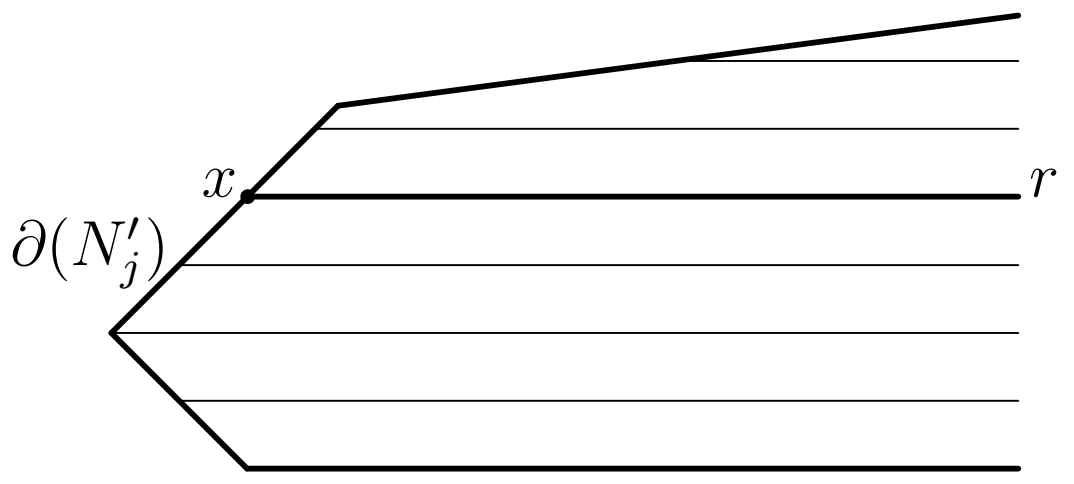

Bad case: no continuous retraction since points of $\partial(N'_j)$ are not to move

Nevertheless, modifying the lines of motion slightly such that points are departing from $r$ when they are approaching $\partial(N'_j)$ works.

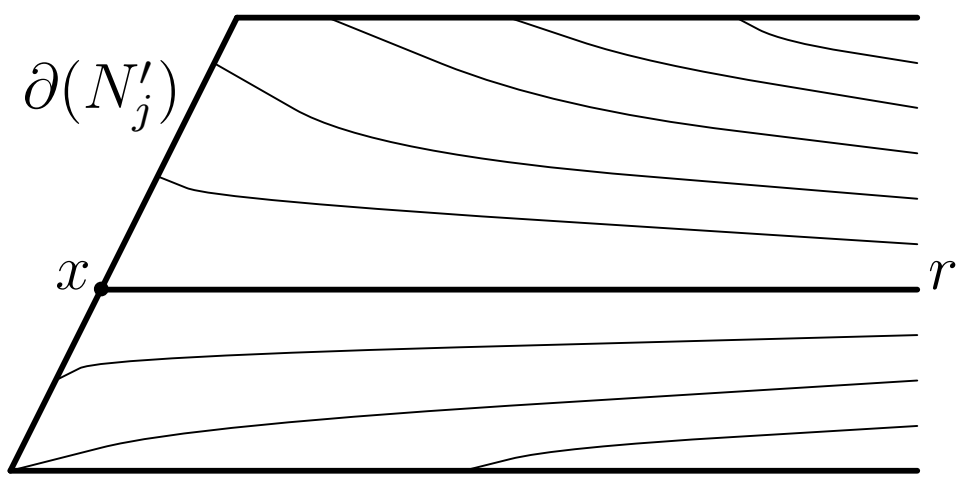

Modified lines of motion yielding a retraction

We just have to make sure that these lines intersect the boundary $\partial(N_j')$. This is achieved most easily by ensuring that they even intersect a supporting hyperplane of $N_j'$ at $x$ which can be accomplished, e.g., by choosing the lines of motion to form a family of affine hyperbolas.

$N_j'$ is bounded: Topologically $N_j'$ is, a disc. Since we are assuming that it has interior points, its dimension is the dimension of the surrounding space $\mathbb{E}$ which equals the dimension of its dual $\mathbb{E}^*$.

In this case, furthermore, 0 lies in the convex hull of $\Lambda'$. By Caratheodory's Theorem [Eggl58, Remark (ii), page 38], there are $\dim \mathbb{E}^* + 1$ elements in $\Lambda'$ whose convex hull contains 0. Since $\Lambda$ and therefore $\Lambda'$ are $m$-tame, $m < \dim \mathbb{E}^* + 1$ follows. Therefore $N_j'$ is a disc of dimension at least $m$.

In either case, attaching $N_j'$ does not introduce nontrivial elements in the homotopy groups up to $\pi_{m-1}$.

Suppose now $\Lambda'$ does not span $\mathbb{E}^*$. Then we write $N_j'$ as a direct product of

$$V := \bigcap_{\lambda \in \Lambda'} \ker(\lambda)$$

and a convex subset $C \subseteq \mathbb{E}/V$. This quotient is dual to the span of $\Lambda'$, and the considerations of the previous paragraph apply to $C$, yielding the same result. **q.e.d.**

In order to apply this proposition to $\boldsymbol{Y}$ we have to prove

**Lemma 4.12.** *The set* $\left(\biguplus_{v \in S} \Phi_v^+\right)\big|_H$ *is* $(|S|-1)$*-tame.*

This set is certainly not $|S|$-tame because of the product formula. Thus the lemma roughly says that there are no further relations among the places.

**Proof.** Let

$$\sum_{v \in S, \alpha \in \Phi^+} \mu_{v,\alpha} \alpha_v \tag{1}$$

be a convex combination, that is, all coefficients $\mu_{v,\alpha}$ are nonnegative and add up to 1. We have to show that this combination does not vanish on $H$ unless at least $|S|$ coefficients are $\neq 0$. In fact, we will see that for each prime $v \in S$ there is at least one non-vanishing coefficient $\mu_{v,\alpha}$.

Call a prime $v$ *void* if all coefficients $\mu_{v,\alpha}$ indexed by that prime vanish. Any positive root is a positive integer combination of base roots. Substitution of all the positive roots by these combinations in (1) yields again a nontrivial positive combination vanishing on $H$, now involving only base roots. Note that this process may increase or decrease the number of vanishing coefficients, but it does not alter the set of void primes.

Therefore, since we want to prove that the set of void primes is empty, we may assume without loss of generality that we are dealing with a convex combination

$$\sum_{v \in S, \alpha \in \Delta} \mu_{v,\alpha} \alpha_v \tag{2}$$

of base roots. If this combination vanishes on $H$ we conclude that

$$\sum_{v \in S, \alpha \in \Delta} \mu_{v,\alpha} \alpha_v \in \operatorname{im}(\boldsymbol{\zeta}^*) .$$

Hence, for each base root $\alpha$ involved, we have

$$\mu_{v,\alpha} = \mu_{v',\alpha} \; \forall v, v' \in S$$

by definition of $\boldsymbol{\zeta}$. Hence none of these coefficients vanishes because we started with a nontrivial combination. **q.e.d.**

## 4.3 Negative Results

In order to complete the proof of Theorem B we have to show that $\mathcal{B}(\mathcal{O}_S)$ is not of type $\mathrm{FP}_{|S|}$. We proof this by reduction to the rank-1-case, but we present only the algebraic reduction that applies to the Borel subgroups $\mathrm{B}_n^0 \leq \mathrm{SL}_n$ and $\mathrm{B}_n \leq \mathrm{GL}_n$. The general case is dealt with by an analogous geometric reduction to the rank-1-case [Bux97b, Section 5].

**Theorem 4.13 ([Bux97a, Theorem 8.1]).** *The group $\mathrm{B}_2^0(\mathcal{O}_S)$ is not of type $FP_{|S|}$.*

This can be proven by various means. Since $\mathrm{B}_2^0(\mathcal{O}_S)$ is metabelian, one can use the $\Sigma$-theory of geometric invariants introduced by R. Bieri and R. Strebel [BiSt80]. This line of reasoning was taken in [Bux97a].

Another way is via Bestvina-Brady-Morse theory [BeBr97]. For $\mathrm{B}_2^0(\mathcal{O}_S)$ the Euclidean buildings $X_v$ are trees and the map

$$\boldsymbol{\pi} : \boldsymbol{X} \to \mathbb{R}^r = \mathbb{R}$$

takes values in the real line since $\mathrm{SL}_2$ has rank $r = 1$. So we regard it as a *height* on $\boldsymbol{X}$. To this map, Bestvina-Brady-Morse theory applies, which can be used to reprove that $\mathrm{B}_2^0(\mathcal{O}_S)$ is of type $\mathrm{F}_{|S|-1}$ as follows. Ascending and descending links are easily computed because of the product structure of $\boldsymbol{X}$. They turn out to be joins of ascending, respectively descending, links in the factors. Hence they are points or joins of $|S|$-spheres, and in either case $(|S|-1)$-connected. From this the claim follows by the Morse lemma of [BeBr97].

The method can be refined to yield as well a

**Proof of Theorem 4.13.** Consider the family of compact intervals of $\mathbb{R} = \mathbb{R}^r$ as a directed system via the relation of inclusion. We claim that the corresponding directed system

$$(\pi_{|S|-1}(\boldsymbol{\pi}^{-1}(I)))_{I\subset\mathbb{R}}$$

of homotopy groups of preimages in $\boldsymbol{X}$ is not essentially trivial.

To show that, for any $I \ni 0$, the map

$$\pi_{|S|-1}(\boldsymbol{\pi}^{-1}(0)) \to \pi_{|S|-1}(\boldsymbol{\pi}^{-1}(I))$$

is not trivial, we take a vertex $\boldsymbol{x} \in \boldsymbol{X}$ above $\max(I)$. In its descending link there are spheres that are nontrivial in the whole link of $\boldsymbol{x}$. Pick one and project it along geodesics until it is on height 0. We see, that this height-0-sphere is nontrivial in $\boldsymbol{X} - \{\boldsymbol{x}\}$ because $\boldsymbol{X}$ is CAT(0) and there is a geodesic projection $(\boldsymbol{X} - \{\boldsymbol{x}\}) \to \mathrm{Lk}_{\boldsymbol{X}}(\boldsymbol{x})$. This projection takes back the height-0-sphere to the nontrivial homotopy element in $\mathrm{Lk}_{\boldsymbol{X}}(\boldsymbol{x})$ we stated with. Hence, the height-0-sphere does not become trivial in $\boldsymbol{\pi}^{-1}(I)$. **q.e.d.**

**Lemma 4.14.** $\mathrm{B}_n(\mathcal{O}_S)$ *and* $\mathrm{B}_n^0(\mathcal{O}_S)$ *have identical finiteness properties.*

**Proof.** Let $\mathrm{D}_n$ denote the scheme of diagonal matrices of rank $n$ and let $\mathrm{D}_n^0$ denote the subscheme of those matrices that have determinant 1. Note that $\mathrm{D}_n(\mathcal{O}_S)$ is finitely generated and contains a free Abelian subgroup $A$ of finite index that decomposes as a product $A = A^0 \times A'$ where $A^0$ is a free Abelian subgroup of finite index in $\mathrm{D}_n^0(\mathcal{O}_S)$. Moreover, we can chose $A'$ to be a group of multiples of the identity matrix. Hence $A'$ acts trivially on the unipotent radical.

Writing

$$\mathrm{B}_n(\mathcal{O}_S) = \mathrm{D}_n(\mathcal{O}_S) \ltimes \mathrm{U}_n(\mathcal{O}_S)$$

and

$$\mathrm{B}_n^0(\mathcal{O}_S) = \mathrm{D}_n^0(\mathcal{O}_S) \ltimes \mathrm{U}_n(\mathcal{O}_S)$$

we see that $\mathrm{B}_n(\mathcal{O}_S)$ is commensurable to $A \ltimes \mathrm{U}_n(\mathcal{O}_S) = A' \times (A^0 \ltimes \mathrm{U}_n(\mathcal{O}_S))$, which has the same finiteness properties as $A^0 \ltimes \mathrm{U}_n(\mathcal{O}_S)$. This, in turn, is commensurable to $\mathrm{B}_n^0(\mathcal{O}_S)$. **q.e.d.**

Hence there is no loss in confining ourselves to $\mathrm{B}_n(\mathcal{O}_S)$. For this group, the following observation easily accomplishes the the reduction to $\mathrm{B}_2(\mathcal{O}_S)$.

**Observation 4.15.** $\mathrm{B}_2(\mathcal{O}_S)$ *is a retract of* $\mathrm{B}_n(\mathcal{O}_S)$ *in the following way:*

$$\begin{pmatrix} * & * & \cdots & * \\ 0 & * & \ddots & \vdots \\ \vdots & \ddots & \ddots & * \\ 0 & 0 & 0 & * \end{pmatrix}$$

*with the obvious inclusion and projection map.* **q.e.d.**

By Remark 3.3 it follows that $\mathrm{B}_2(\mathcal{O}_S)$ has all finiteness properties that $\mathrm{B}_n(\mathcal{O}_S)$ enjoys. Hence $\mathrm{B}_n(\mathcal{O}_S)$ is not of type $\mathrm{FP}_{|S|}$ in view of Theorem 4.13 whence $\mathrm{B}_n^0(\mathcal{O}_S)$ is neither.

# References

[Abel79] *H. Abels:* An example of a finitely presented solvable group. In: [Wall79] p. 205 – 211 [cited on page 9]

[Abel87] *H. Abels:* Finite presentability of $S$-arithmetic groups – Compact presentability of solvable groups. Springer LNM 1261 (1987) [cited on pages 8, 10 and 14]

[Abel91] *H. Abels:* Finiteness Properties of certain arithmetic groups in the function field case. Israel J. Math. 76 (1991) p. 113 – 128 [cited on page 9]

[AbBr87] *H. Abels, K.S. Brown:* Finiteness Properties of Solvable $S$-Arithmetic Groups: An Example. Journal of Pure and Applied Algebra 44 (1987) p. 77 – 83 [cited on page 9]

[AbTi97] *H. Abels, A. Tiemeyer:* Compactness Properties of Locally Compact Groups. Transformation Groups 2 (1997) p. 119 - 135 (confer [Tiem97]) [cited on pages 10, 11, 12 and 27]

[Abra87] *P. Abramenko:* Endlichkeitseigenschaften der Gruppen $\mathrm{SL}_n(\mathbb{F}_q[t])$. Thesis Frankfurt am Main (1987) [cited on page 9]

[Abra96] *P. Abramenko:* Twin Buildings and Applications to $S$-Arithmetic Groups. Springer LNM 1641 (1996) [cited on page 9]

[Behr62] *H. Behr:* Über die endliche Definierbarkeit von Gruppen. Journal für die reine und angewandte Mathematik (Crelles Journal) 211 (1962) p. 116 – 122 [cited on page 8]

[Behr98] *H. Behr:* Arithmetic groups over function fields I. A complete characterization of finitely generated and finitely presented arithmetic subgroups of reductive algebraic groups. Journal für die reine und angewandte Mathematik (Crelles Journal) 495 (1998) p. 79–118 [cited on page 9]

[Behr99] *H. Behr:* Higher finiteness properties of $S$-arithmetic groups in the function field case I: Chevalley groups with coefficients in $\mathbb{F}_q[t]$. Preprint, Frankfurt (1999) [cited on page 9]

[BeBr97] *M. Bestvina, N. Brady:* Morse theory and finiteness properties of groups. Inventiones mathematicae 129 (1997) p. 445 – 470 [cited on pages 4 and 23]

[BiSt80] *R. Bieri, R. Strebel:* Valuations and Finitely Presented Metabelian Groups. Proceedings of the London Mathematical Society 41 (1980) p. 439 – 464 [cited on pages 9 and 23]

[Bore63] *A. Borel:* Some finiteness properties of adele groups over number fields. Publ. Math., Inst. Hautes Etud. Sci. 16 (1963) p. 5 – 30 [cited on page 13]

[Bore91] *A. Borel:* Linear Algebraic Groups. Second Enlarged Edition; Springer GTM 126, New York (1991) [cited on page 7]

[BoHa62] *A. Borel, Harish-Chandra:* Arithmetic subgroups of algebraic groups. Annals of Mathematics 75 (1962) p. 485 – 535 [cited on page 8]

[BoSe73] *A. Borel, J-P. Serre:* Corners and arithmetic groups. Commentarii Mathematici Helvetici 48 (1973) p. 436 – 491 [cited on page 8]

[BoSe76] *A. Borel, J.-P. Serre:* Cohomologie d'immeubles et de groupes S-arithmétiques. Topology 15 (1976) p. 211 – 232 [cited on pages 8, 14 and 15]

[BoTi65] *A. Borel, J. Tits:* Groupes réductifs. Publ. Math., Inst. Hautes Etud. Sci. 27 (1965) p. 55 – 150 [cited on pages 14 and 15]

[Brow82] *K.S. Brown:* Cohomology of Groups. Springer GTM 87, New York (1982) [cited on page 5]

[Brow87] *K.S. Brown:* Finiteness Properties of Groups. Journal of Pure and Applied Algebra 44 (1987) p. 45 – 75 [cited on pages 4 and 11]

[Brow89] *K.S. Brown:* Buildings. Springer, New York, Berlin (1989) [cited on page 16]

[BrTi72] *F. Bruhat, J. Tits:* Groupes Réductives sur un Corps Local I. Publications Mathématiques IHES 41 (1972) p. 5 – 252 [cited on page 16]

[Bux97a] *K.-U. Bux:* Finiteness Properties of some Metabelian $S$-Arithmetic Groups. Proceedings of the London Mathematical Society 75 (1997) p. 308 – 322 [cited on pages 9, 20 and 23]

[Bux97b] *K.-U. Bux:* Endlichkeitseigenschaften auflösbarer arithmetischer Gruppen über Funktionenkörpern. PhD thesis (Frankfurt, 1997) http://math.kubux.net/math/bux98.ps.gz [cited on pages 9, 10, 17, 18, 20 and 23]

[Bux99] *K-U. Bux:* Finiteness Properties of Soluble Arithmetic Groups over Global Function Fields. Preprint (1999) http://math.kubux.net/math/bux99b.ps.gz [cited on pages 9 and 10]

[CaFr67] *J.W.S. Cassels, A. Fröhlich:* Algebraic Number Theory. Proceedings of an Instructional Conference organized by the London Mathematical Society (A Nato Advanced Study Institute) with Support of the International Mathematical Union. Academic Press, London (1967) [cited on pages 10 and 17]

[Chev60] *C. Chevalley:* Certain Schémas de Groupes Semi-Simples. Séminaire Bourbaki 13e année (1960/61) n∘ 219 p. 1 – 16 [cited on page 16]

[Eggl58] *H.G. Eggleston:* Convexity. Cambridge Tracts in Mathematics and Mathematical Physics, Cambridge (1958) [cited on page 22]

[Heyd99] *A. von Heydebreck:* Homotopy properties of certain complexes associated to spherical buildings. To appear: Israel Journal of Mathematics [cited on page 19]

[Knes64] *M. Kneser:* Erzeugende und Relationen verallgemeinerter Einheitengruppen. Journal für die Reine und Angewandte Mathematik 214/215 (1964) p. 345 – 349 [cited on page 8]

[Macb64] *A.M. Macbeath:* Groups of homeomorphisms of a simply connected space. Annals of Math 79 (1964) p. 473–488 [cited on page 1]

[O'Me73] *O.T. O'Meara:* Introduction to Quadratic Forms. Berlin, Heidelberg, New York (Springer, 1973) [cited on pages 10 and 17]

[Quil78] *D. Quillen:* Homotopy Properties of the Poset of Nontrivial $p$-Subgroups of a Group. Advances in Mathematics 28 (1978) p. 101 – 128 [cited on page 4]

[Ragh68] *M.S. Raghunatan:* A Note on Quotients of Real Algebraic Groups by Arithmetic Subgroups. Inventiones mathematicae 4 (1968) p. 318 – 335 [cited on page 8]

[Rona89] *M. Ronan:* Lectures on Buildings. Perspectives in Mathematics, Academic Press, San Diego, London (1989) [cited on page 16]

[Serr77] *J.-P. Serre:* Arbres, amalgames, $SL_2$. Asterique 46 (1977) [cited on page 1]

[Span66] *E.H. Spanier:* Algebraic Topology. Springer, New York (1966) [cited on page 5]

[Stei68] *R. Steinberg:* Lectures on Chevalley Groups, Notes prepared by John Faulkner and Robert Wilson.. Yale University, New Haven, Connecticut, (1968) [cited on page 16]

[Stuh76] *U. Stuhler:* Zur Frage der endlichen Präsentiertheit gewisser arithmetischer Gruppen im Funktionenkörperfall. Mathematische Annalen 224 (1976) p. 217 – 232 [cited on page 9]

[Stuh80] *U. Stuhler:* Homological properties of certain arithmetic groups in the function field case. Inventiones mathatematicae 57 (1980) p. 263 – 281 [cited on page 9]

[Tiem97] *A. Tiemeyer:* A Local-Global Principle for Finiteness Properties of $S$-Artihmetic Groups over Number Fields. Transformation Groups 2 (1997) p. 215 – 223 (confere [AbTi97]) [cited on pages 10, 13, 14 and 25]

[Wall79] *C.T.C. Wall:* Homological Group Theory. LMS Lecture Notes 35, Cambridge University Press (Cambridge 1979) [cited on page 24]

[Weil73] *A. Weil:* Basic Number Theory. Reprint of the Second Edition (1973); Springer Classics in Mathematics, Berlin Heidelberg (1995) [cited on pages 10 and 17]

I would like to thank Prof. Heinz Helling and Prof. Thomas Müller for giving me the opportunity to present part of this work at the conference they organized at Bielefeld in August 1999. I am indebted to Peter Brinkmann and Blake Thornton for very helpful criticism on the manuscript.

Kai-Uwe Bux
Department of Mathematics
Cornell University
310 Malott Hall
Ithaca, NY 14853-4201

email: `bux_math@kubux.net`